# Collatz Conjecture: Exposition and Proof through a Structured Approach

Ken Surendran*     and     Desarazu Krishna Babu**

**Abstract**

A structured approach for the Collatz conjecture is presented using just the *odd* integers that are, in turn, divided into categories based on the roles they play such as Starter, Intermediary and Terminal. The expression *4x+1* is used as a tool to expose all the hidden and significant characteristics of the conjecture that lead us to its proof. The *mixing* properties of the *iterates* are addressed by showing that the *Collatz iterates* of half of <u>all</u> the *odd* integers that are of the form 4*m*+3, on the average, increase by three times the value of the *odd* integer that was used to start with, while the *iterates* of those of 4*m*+1, on the average, decrease by a factor of four. Further, expressions are provided to generate all the sets of *odd* integers where the Collatz *iterate* of all the integers in each set is an integer of the form *6m+1* or *6m+5*. The significance of the Collatz *net (tree)* is obvious since it encompasses all the Collatz *trajectories*.

## 1. Introduction

Let *x* be a positive integer. As per the Collatz conjecture, whenever *x* is *odd,* we replace it with (3*x*+1)/2 and whenever *x* is *even* we replace it with x/2 and ultimately the value of *x* becomes 1 at some stage. So, when *x* is *odd*, compute $y = (3x+1)/2^{\alpha}$, where α is a positive integer and it is chosen such that *y* becomes an *odd* integer. If the resulting *y* is > 1 make *x* = *y* and repeat the above process. As per the conjecture, the value of *y* ultimately reaches 1 in a finite number of steps. This is the essence of the Collatz conjecture, also known as the *3x+1 Problem*. The initial exposure to the Collatz conjecture topic was through Wikipedia [1]. Later, the book [2] edited by J. C. Lagarias motivated us to work on this *3x+1 Problem* and we found the book to be a comprehensive resource on this topic. In this study, we consider all the *odd* positive integers only. Also here, 'conjecture' refers to the *Collatz conjecture*.

In Section 2, we describe, from the conjecture's perspective, a few categories of the *odd* integers. In one categorization, we use the same type of classification of *odd* integers (albeit, with different names) that is discussed by Tao in his recent publication [3]. This helps us to come up with a new type of *Collatz graph* with a specific set of (i) *Starter* integers, with which the conjecture process begins and (ii) the *Intermediary* integers, which form the core and the essence of the process. The *Intermediary integers* are the *odd iterates* (*i.e.*, the conjecture results) found in all the *trajectories*. In addition, we also introduce the *Terminal* integers, which end the process by producing 1 as their conjecture result.

In Section 3, we start with sample conjecture *trajectories*, presented in a Table, for a few *Starter* integers and also in the form of graphs just with the *odd iterates*. In the graphs, we observe the presence of 4*x*+1 that links the *iterates* from the conjecture *trajectories* of the different *Starter* integers.

___________________________
* Southeast Missouri State Univ., Dept. of Computer Science (Retired; *suren@)linuxmail.org*)
** Mobil Oil Corporation, Senior Engineering Advisor (Retired; earlier with Princeton Univ.)

Based on this, in Sections 4, we develop a generalized version of 4*x*+1 and use it in Section 5 with *all* the *odd* integers to explore the hidden characteristics of the conjecture itself. We provide sample Tables where each row has a set of *odd* integers, all of which produce the same conjecture *iterate* and, most importantly, we provide expressions to get such sets of *odd* integers for **all** of the *odd iterates*. That is, each row in these *extended* Tables has a *predecessor set* (see page 62 in [2]) of the corresponding *Intermediary* integer in the conjecture *iterate* column. We also provide the basic logic for constructing the *conjecture trajectory* for any *odd* integer in a different way by accessing the results in these Tables using suitable expressions.

In Section 6, we also introduce the *Collatz Tree* with its root (*i.e.*, *End* integer 1) above, the branches (*i.e.*, the *Intermediary* integers) below, and the leaves (*i.e.*, *Starter* integers) spread out everywhere. As a result, **all** the *odd* integers go up towards the root using the conjecture. Here, we start with known integers for the first two layers below the root (*i.e.*, 1) of the *Collatz Tree* and provide the core logical steps for constructing the remaining lower layers with several segments (one under each *iterate* in the layer above).

In Section 7, we study the complete behavior of all the *odd* integers that use just one division by 2 in their conjecture process that results in sequentially increasing *iterates*. Half of all the *odd* integers are in this category. We classify, using an expression, each of these integers based on its *total* number of sequentially increasing *Collatz iterates* and present all the results in a Table. (We use the term *Alpha1 Trajectory Length* for such a *total* number; *Alpha1* is: $\alpha = 1$).

In Section 8, based on such detailed studies with all the *odd* integers and their *iterates*, we show with realistic (not *heuristic*) evidence that, for half of **all** the *odd* integers, their *odd iterates* in the corresponding *trajectories*, on the average, increase (swell) by a factor of 3 while for the other half of **all** the *odd* integers the corresponding *trajectories* reduce (shrink), on the average, by a factor of 4. This clarifies the *mixing* properties and proves that there cannot be any divergent *trajectories*. In Section 9, following a brief summary of this study, we emphasize the significance of the Collatz net (*tree*) since it encompasses all the strings (*trajectories*) of the net.

## 2. Types of odd Integers from Collatz Conjecture Perspective

We start with 1 since that is where everything *ends* in the Collatz conjecture. Applying the conjecture to 1: *i.e.,* 3 X 1 + 1 = 4; 4/2 =2; and 2/2= 1, we note that the process keeps looping. To progress further, we add the result before the division (*i.e.* 4) to the starting integer (*i.e.,* 1) and apply the conjecture to this sum (1+4 = 5): 3 X 5+1 =16 = $4^2$ and $16/2^4 = 1$. Continuing with the same approach (*i.e.*, adding the result before the division to its starting integer), the new sum is: 5+16 = 21; and 3X21+1= 64 = $4^3$ and $64/2^6 = 1$. Generalizing the above approach, let $T_k$, (with $k = 0, 1, 2, 3,$ ) be a set of integers generated by:

$$T_k = 1+4+4^2+4^3 + \ldots+4^k = \sum_{n=0}^{k}(4^n) = (4^{k+1} - 1)/3 \qquad (1)$$

Applying the conjecture to $T_k$, the final result is 1 as summarized below:

$$3T_k + 1 = 4^{k+1} \text{ and } 4^{k+1}/2^{2(k+1)} = 1 \qquad (2)$$

The integers $T_k s$ (for every positive integer *k*) collectively provide the required closure to the conjecture. We call them the *Terminal* integers (consisting of 1, 5, 21, 85, 341, 1365, …, etc.) since they terminate the Collatz conjecture process by producing 1, the *End* integer, as their result. These *Terminal* integers are also generated, using a well-known method [1] (*aka* a Syracuse Function property) starting with $T_0 = 1$ and using

$$T_{k+1} = 4T_k + 1. \qquad (3)$$

Now, let us apply the conjecture for a few *odd* integers: 3, 9, 15 and 21. Their *trajectories* with only the *odd iterates* [using $(3x+1)/2^\alpha$] are: **3** → 5 → 1; **9** → 7 → 11 → 17 → 13 → 5 → 1; **15** → 23 → 35 → 53 → 5 → 1; and **21** → 1. Both 5 and 21, being *Terminal* integers, end the *trajectory* with 1 in one step. Both 13 and 53, as they go through the conjecture, result in 5. The rationale for this lies in the generalized version of Equation (3) [53= 4 X 13+1]:

$$y = 4x+1 \qquad (4)$$

Equation (4) will be used extensively later on. Since the $3x+1$ conjecture is about 3, we first examine – using a Theorem – the significance of the *odd* integers that are '*odd* multiples of 3'.

THEOREM 2.1. In the conjecture process, an *odd iterate* cannot be an *odd* multiple of 3.
*Proof.* Let $x$ and $y$ be *odd* integers such that $y = (3x+1)/2^\alpha$, where α, a positive integer, is chosen such that $y$ is an *odd* integer. We reverse the process and check if $y$ could be an *odd* multiple of 3 (*i.e.,* $y = 3z$, where $z$ is an *odd* integer: 1, 3, 5, …) and examine what type of $x$ we get as a result of this new restriction on $y$:

$$x = (2^\alpha . X\ y - 1)/3 = (2^\alpha\ X\ 3z - 1)/3 = 2^\alpha z - 1/3 \qquad (5)$$

We notice that, for any *odd* value of $z$ and for any positive integer value of $\alpha$, $x$ could no longer be an integer due to the restriction we placed on $y$. Hence we realize that, for any *odd* integer values of $x$, the *result* of $(3x+1)/2^\alpha$ can never be an *odd* multiple of 3. So, the *odd* multiples of 3 could only be at the <u>starting</u> positions of the Collatz conjecture process     (*Q.E.D.*)

**2.1. Starter and Intermediary Integers.** In view of *Theorem-2.1*, the *odd* integers are split into two types: (i) the *odd* multiples of 3: $3(2m + 1)$ for $m = 0, 1, 2, 3, $ , (we name them as the *Starter* integers) and (ii) all the rest of the *odd* integers: $3(2m + 1)\pm 2$, with $m = 0, 1, 2, 3, $ ,: we name them as the *Intermediary* integers and they are of two types: *6m+1* and *6m+5*. In the Collatz conjecture, we begin the process with each of the *Start*er integers and continue applying the conjecture sequentially to the resulting *odd* integers until we reach a *Terminal* integer which then leads to the *End* integer 1. All those resulting *odd* integers occurring in between the *Starter* and *End* integers (including the *Terminal* integer) are found to be the *Intermediary* integers. [Here is an example to illustrate the classification of *odd* integers from the Collatz conjecture perspective: Earlier on, we saw the Collatz conjecture *trajectory* for integer 9: **9** → 7 → 11 → 17 → 13 → 5 → 1. Here, 9 is a *Starter* integer, 5 is a *Terminal* integer and 1 is the *End* integer. The other in-between (*i.e.,* 7, 11, 17, 13) are the *Intermediary* integers. Here 7 and 13 are of type *6m+1;* and 11, and 17 are of type *6m+5*] We now introduce the idea of a *Conjecture Table* using a set of *Terminal* integers for illustration.

**2.2. Collatz Conjecture Table for Terminal Integers.** We present in Table-1 with the columns showing the *iterates* of the conjecture applied to *odd* integers beginning with the corresponding *Starter* integers for a set of initial *Terminal* integers (5, 21, 85, 341, 1365,…). We note that, starting with 21, every third *Terminal* integer is also a *Starter* integer (*odd* multiple of 3). In Table-1, the *Terminal* integers are underlined and they all converge to 1 directly with no in-between *iterates*.

Using a *reverse process* – i.e., going up the *trajectory* – we can find the *Starter* integer for a given *Intermediary-cum-Terminal* integer like 85. We consider, for n ≥ 1, $2^n$X85 and subtract 1: *i.e.,* 2X85 –1 = 169, 4X85–1 = 339, 8X85–1= 679; and pick the first result that is a multiple of 3. Here, we see that 339 (= 3X113) is a multiple of 3. However, we see that 113 is not a *Starter* integer and so we repeat the process now with 113: (2X113 – 1) = 225 (= 3X75). Since 75 is a *Starter* integer, the *reverse* process ends here.

| Starter | 3 | 21 | 75 | 201 | 1365 | 7281 | 17019 | *87381* | 245481 | 932067 | 5592405 |
|---|---|---|---|---|---|---|---|---|---|---|---|
| Intermediary | | | 113 | 151 | | | 25529 | | 184111 | | |
| Intermediary | | | | 227 | | | 19147 | | 276167 | | |
| Intermediary | | | | | | | 14563 | | 414251 | | |
| Intermediary | | | | | | | | | 621377 | | |
| Intermediary | | | | | | | | | 466033 | | |
| Terminal | 5 | 21 | 85 | 341 | 1365 | 5461 | 21845 | 87381 | 349525 | 1398101 | 5592405 |
| End | 1 | 1 | 1 | 1 | 1 | 1 | 1 | 1 | 1 | 1 | 1 |

TABLE 1. Collatz Conjecture Table with Sample Terminal Integers with Their Starter Integers

In the following, we present a Table that shows the complete *Collatz trajectory* for a few smaller *Starter* integers to understand the behavior of the *iterates*. We also present the *Collatz graphs*, drawn using the results in this Table. The objective is to identify, through observation and analysis of these results, any hidden intrinsic characteristics of the Collatz conjecture.

## 3. Collatz Conjecture Tables and Graphs Beginning with A Few Starter Integers

Having defined, from the conjecture perspective, the three main types of integers (*Starter, Intermediary and Terminal*), we now apply the conjecture to a set of sample *Starter* integers (from 3 to 255). We show the results consisting of only the *odd* integers in Table 2. The first column shows the *Starter* integer from 3 to 255; here, some of the *odd* integers with shorter Collatz *trajectories* are presented in the same row. All the *Starter* integers, as they go through the conjecture, end with one of the *Terminal* integers followed by the *End* integer 1. The *Intermediary* integers, which are the *iterates* in the *trajectory*, are listed in the second column in the order in which they occur. For convenience, we use 27 as the reference since some *Starter* integers (> 27) have several common *Intermediary* integers that are in 27's conjecture *trajectory*. In the following we present a few observations from Table 2:

 (i) In several of these *Collatz trajectories*, there are continuous steady increase in the values of conjecture *iterates* by about 1½ times the previous value, since they all use (3x+1) /2 *i.e.*, with $\alpha$ = 1. For instance, in the *Collatz trajectory* for the *Starter* integer 255, there are *seven* continuous increases before a reduction occurs:

**255** → 383 → 575 → 863 → 1295 → 1943 → 2915 → 4373 → 205→ … → 1

However, in this case, this steady increase is followed by a drastic reduction with $\alpha$ = 6, taking the *iterate* value lower than the one we started with (*i.e.*, 255 to 205). It appears that there is a need to study the number of such continuous *iterate* increases for all such *odd* integers..

(ii) The *Intermediary* integers occur only once in each *trajectory*

(iii) The lengths of the *Collatz trajectories* vary considerably.

(iv) The *Starter* integers use different values of α in computing their first Collatz conjecture *iterate*. We note that α = 1 is used for every alternate *odd* multiples of 3 starting from 3 (*i.e.,* 3, 15, 27, …); α = 2 for every fourth starting from 9 (*i.e.,* 9, 33, 57, …); and so on. Next, we present a new type of *Collatz Conjecture graph* using the results in Table 2.

**3.1. Collatz Conjecture Graphs.** Based on the results in the Table 2, we present a new type of *Collatz Graphs* (different to those presented in [1] and [2]) and we use *MS-Visio* to draw these graphs. In these flow-chart form of graphs, we use different graphic symbols to place the different types of *odd* integers. In our Tables, there are just 43 *Starter* (3 to 255) integers. For the sake of convenience in presenting the conjecture graphs, we split them into three groups. First, we collected all those *Starter* integers that have 5 as their *Terminal* integer and put them into two categories: (i) those that reach 5 from *Pre-Terminal* 13 (13 → 5) and (ii) the rest that reach 5 from *Pre-Terminal* 53 (53 → 5). (The *Pre-terminals* of 5 are: 3, 13, 53, 213, …→ 5.) Also, we put all those integers that do not have 5 as their *Terminal* integers into category (iii). (Note: In these graphs, the flow is assumed to be downwards, when there are no arrows.)

Figure 1A shows the *Collatz Graph* for categories (i) and (iii) and Figure 1B shows the graph for category (ii). Also, the symbols used for the different types of integers are explained under the graphs. In these graphs, some of the *iterates* are put together in a single block mainly for convenience and also the graphs (i) and (ii) can be combined into one with just a single *Terminal* integer box for 5. In these *Collatz Graphs*, every *odd* integer occurs only once. This is expected since these graphs are based on the results presented in Table 2 where there are no duplicate *iterates* within any specific *Collatz trajectory*. In the following, we discuss some significant observations made from these *Collatz Graphs*.

**3.2 Observations and Presence of 4x+ 1.** We see the role of Equation (4) in the two *Collatz Graphs* (Figures 1A and 1B). For instance in Fig. 1B, the integers 27, 109, and 437 produce the same *iterate* 41 as they go through the conjecture process: (3X27 +1)/2 = (3X109+1)/8 = (3X437+1)/32 = 41. Here we note the presence of Equation (4): 4X27 +1 = 109 and 4X109 +1 = 437. In the two main graphs (Fig. 1A and 1B), we see the presence of Equation (4) at some 34 intersections  In Table 3, we list the two or three integers at these intersections along with their Collatz conjecture *iterate* in the Result-columns.

*Observation*: We considered only 43 *Starter* integers in Table-2. In that small sample, among the 34 *Intermediary* integers in the two Result-columns in Table-3, we note that 11 are of the type 6*m*+1 and 23 are of the type 6*m*+5 (roughly, a ratio of 1:2). Even though 6*m*+1 and 6*m*+5 are equal in number (each with one-third of all the *odd* integers), it appears that relatively more *odd* integers, possibly 2/3$^{rd}$ of them, generate *iterates* of the type 6*m*+5 and the remaining 1/3$^{rd}$ only generate the *iterates* of the type 6*m*+1.

This could be an important observation concerning the *iterates* and the *odd* integers that produce them. Next we examine the importance of 4*x*+1 in this study. To begin with, we look at *odd* integers, expressed in binary, that are related by 4*x*+1.

| Starter Integers | Collatz Conjecture Trajectories [(3x+1)/2^α] for Some Initial Starter Integers | | | |
|---|---|---|---|---|
| 3; 9; 15; 21 | **3** -> 5, 1; | **9**-> 7, 11, 17, 13, 5, 1; | **15**-> 23, 35, 53, 5, 1; | **21**-> 1 |
| 27 | **27**-> 41, 31, 47, 71, 107, 161, 121, 91, 137, 103, 155, 233, 175, 263, 395, 593, 445, 167, 251, 377, 283, 425, 319, 479, 719, 1079, 1619, 2429, 911, 1367, 2051, 3077, 577, 433, 325, 61, 23, 35, 53, 5, 1 | | | |
| 33; 39 | **33** -> 25, 19, 29, 11, 17, 13, 5, 1 ; | | **39** -> 59, 89, 67, 101, 19, 29, 11, 17, 13, 5, 1. | |
| 45; 51; 57 | **45** -> 17, 13, 5, 1; | **51**-> 77, 29, 11, 17, 13, 5, 1; | **57** -> 43, 65, 49, 37, 7, 11, 17, 13, 5, 1. | |
| 63 | **63** - > 95, 143, 215, 323, 485, 91, … **see under 27 for the rest starting in its first line** | | | |
| 69, 75, 81 | **69** -> 13, 5, 1; | **75** -> 113, 85, 1; | **81** -> 61, 23, 35, 53, 5, 1; | |
| 87, 93, 99 | **87** -> 131, 197, 37, 7, 11, 17, 13, 5, 1; | **93** -> 35, 53, 5, 1; | **99** -> 149, 7, 11, 17, 13, 5, 1. | |
| 105 | **105** -> 79, 119, 179, 269, 101, 19, 29, 11, 17, 13, 5, 1; | | | |
| 111 | **111** -> 167, 251, … **see under 27 for the rest starting in its second line** | | | |
| 117, 123 | **117** -> 11, 17, 13, 5, 1; | **123**-> 185, 139, 209, 157, 59, 89, 67, 101, 19, 29, 11, 17, 13, 5, 1. | | |
| 129 | **129** -> 97, 73, 55, 83, 125, 47, … **see under 27 for the rest starting in its first line** | | | |
| 135, 141 | **135** -> 203, 305, 229, 43, 65, 49, 37, 7, 11, 17, 13, 5, 1; | | **141** -> 53, 5, 1. | |
| 147 | **147** - > 221, 83, … **see under 129 and see under 27 for the rest starting in its first line** | | | |
| 153 | **153** -> 115, 173, 65, 49, 37, 7, 11, 17, 13, 5, 1; | | | |
| 159 | **159** -> 239, 359, 539, 809, 607, 911, … **see under 27 for the rest starting in its second line** | | | |
| 165 | **165** -> 31, … **see under 27 for the rest starting in its first line** | | | |
| 171 | **171** -> 257, 193, 145, 109, 41, … **see under 27 for the rest starting in its first line** | | | |
| 177 | **177** ->133, 25, 19, 29, 11, 17, 13, 5, 1. | | | |
| 183 | **183** -> 275, 413, 155, … **see under 27 for the rest starting in its first line** | | | |
| 189 | **189** -> 71, … **see under 27 for the rest starting in its first line** | | | |
| 195 | **195** -> 293, 55, … **see under 129 and see under 27 for the rest starting in its first line** | | | |
| 201 | **201** -> 151, 227, 341, 1. | | | |
| 207 | **207** -> 311, 467, 701, 263, … **see under 27 for the rest starting in its first line** | | | |
| 213, 219 | **213** -> 5, 1. ; | **219** -> 329, 247, 371, 557, 209, 157, 59, 89, 67, 101, 19, 29, 11, 17, 13, 5, 1. | | |
| 225 | **225** -> 169, 127, 191, 287, 431, 647, 971, 1457, 1093, 205, 77, 29, 11, 17, 13, 5, 1. | | | |
| 231 | **231** -> 347, 521, 391, 437, 41, … **see under 27 for the rest starting in its first line** | | | |
| 237 | **237** -> 89, 67, 101, 19, 29, 11, 17, 13, 5, 1. | | | |
| 243 | **243** -> 365, 137, … **see under 27 for the rest starting in its first line** | | | |
| 249 | **249** -> 187, 281, 211, 317, 119, 179, 269, 101, 19, 29, 11, 17, 13, 5, 1. | | | |
| 255 | **255** -> 383, 575, 863, 1295, 1943, 2915, 4373, 205, 77, 29, 11, 17, 13, 5, 1. | | | |

TABLE 2.  Collatz Conjecture Trajectory Table for *Odd* Multiples of 3 from 3 to 255

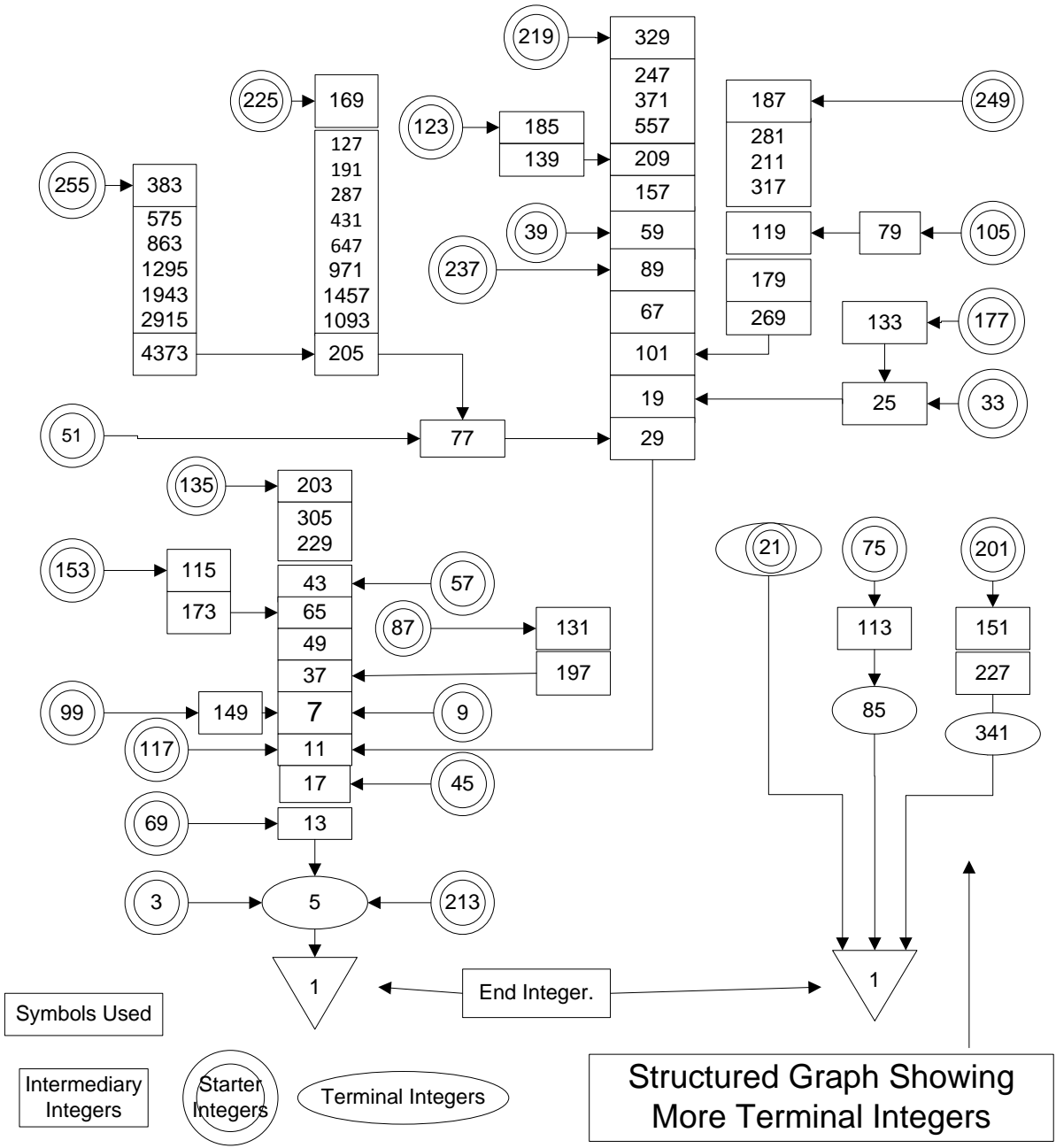

Figure-1A: Example Collatz Structured Graphs

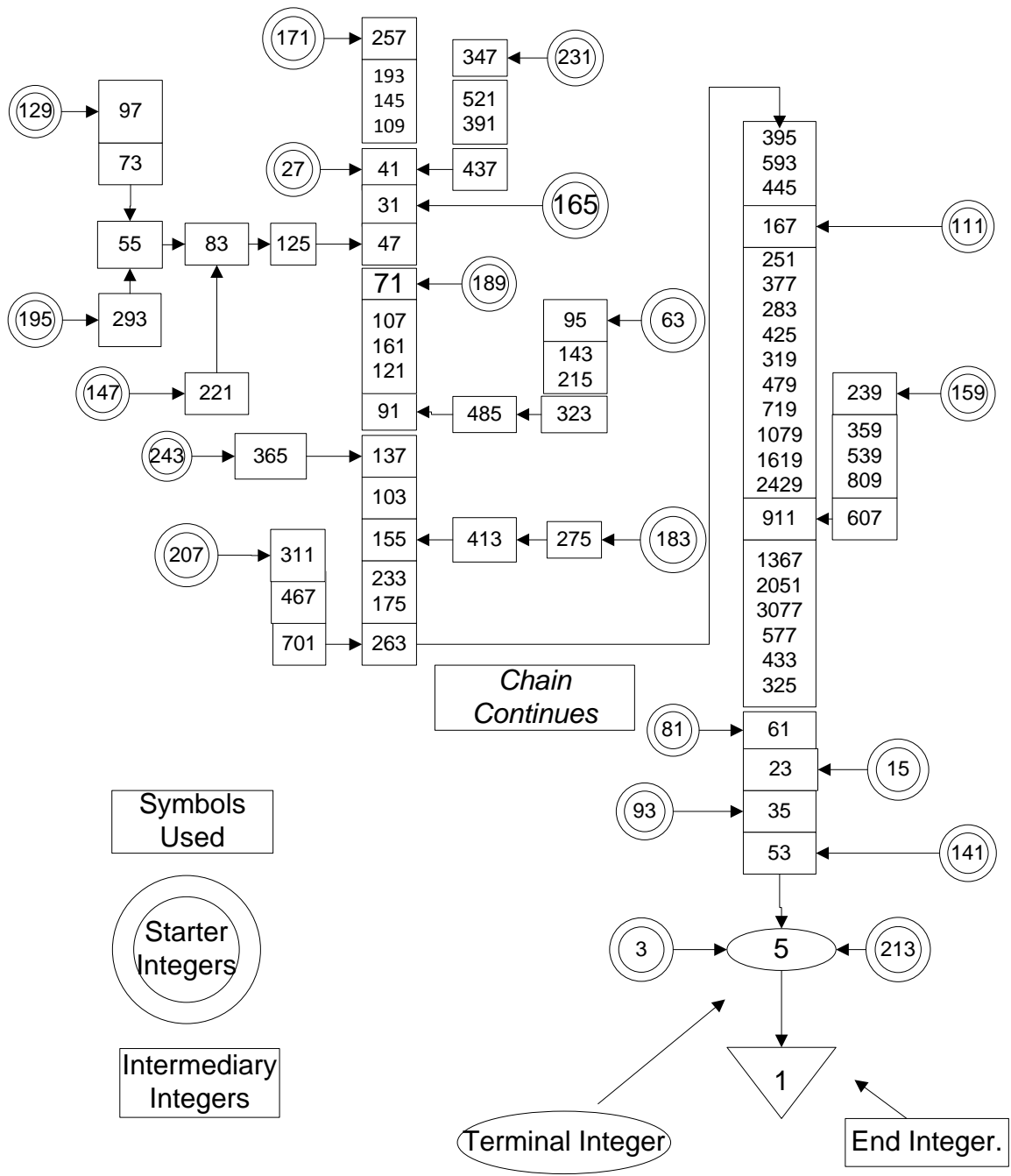

Figure-1B: Example Collatz Structured Graphs

| x | y=4x+1 | z=4y+1 | Result | | x | y=4x+1 | z=4y+1 | Result |
|---|---|---|---|---|---|---|---|---|
| 13 | 53 | 213 | 5 | **** | 39 | 157 | | 59 |
| 9 | 37 | 149 | 7 | **** | 81 | 325 | | 61 |
| 7 | 29 | 117 | 11 | **** | 43 | 173 | | 65 |
| 17 | 69 | | 13 | **** | 47 | 189 | | 71 |
| 11 | 45 | | 17 | **** | 51 | 205 | | 77 |
| 25 | 101 | | 19 | **** | 55 | 221 | | 83 |
| 15 | 61 | | 23 | **** | 59 | 237 | | 89 |
| 33 | 133 | | 25 | **** | 121 | 485 | | 91 |
| 19 | 77 | | 29 | **** | 67 | 269 | | 101 |
| 41 | 165 | | 31 | **** | 79 | 317 | | 119 |
| 23 | 93 | | 35 | **** | 91 | 365 | | 137 |
| 49 | 197 | | 37 | **** | 103 | 413 | | 155 |
| 27 | 109 | 437 | 41 | **** | 111 | 445 | | 167 |
| 57 | 229 | | 43 | **** | 1093 | 4373 | | 205 |
| 31 | 125 | | 47 | **** | 139 | 557 | | 209 |
| 35 | 141 | | 53 | **** | 175 | 701 | | 263 |
| 73 | 293 | | 55 | **** | 607 | 2429 | | 911 |

TABLE 3. Presence of *4x+1* extracted from the two graphs in Fig. 1A and 1B.

### 4. Significance of 4x+1 expression

The *Terminal* integers, $T_k$s, in binary form are: 1, 101, 10101, … (alternate 1 and 0 with 1s at both ends). Multiplying the *Terminal* integers by 3, we get: 11, 1111, 111111,… and adding 1 to this result, we get: 100, 10000, 1000000,…; hence, it is easy to visualize that all $T_k s$ – that are linked by the *4x+1* expression – converge to 1 as they go through the Collatz conjecture. While Equation (1) is intended for *Terminal* Integers, Equations (4) is applicable to all *odd* integers. Let x = 27 (in binary, 11011), its Collatz *iterate* is 41 (in binary 101001). Also *4x +1* = 109 (in binary, 1101101); as 109 goes through the conjecture, it becomes, in binary, 10100100, resulting in 101001 (due to division by $2^2$) which is 41. So, in *4x+1,* the binary form of *x* is preserved and just shifted by two positions to the left.

In the case of *3x+1*, we saw a similar pattern with *Terminal* integers (due to their unique binary format). However, when we compute 3x+1 for non-*Terminal odd* integers, the original binary representation of *x* is distorted. This unique characteristics of 4x+1 leads us to Theorem 4.1 for generalizing an earlier observation in Section 2 (*the odd integers in every trajectory occur only once*) so as to state that there are no duplicate *iterate* within any *Collatz trajectory*.

THEOREM 4.1. As the conjecture process begins with a *non-Terminal Starter* integer and reaches its *Terminal* integer; no duplicates are found among the resulting *Collatz iterates* in the *Collatz trajectory*. (This implies that there are no loops in the *trajectories* as well.)
*Proof*: Let *S* be a *Starter* integer and *T* be its *Terminal* integer in its *Collatz trajectory*: $S \to a \to b \to c \to d \to e \to f \to g \to h \to ..., \to m \to n \to o \to p \to q \to r \to s \to ..., \to T \to 1$. The entries in between *S* and *T* are a set of *Collatz iterates* (*i.e., Intermediary* integers). To prove Theorem 4.1 by negation, we assume that there is a duplication (first one) among these *iterates* say, *g* = *q*. This would mean that, by Equation (4) and by assuming that *p>f*, $p = 4*f + 1$ so that both *f* and *p*, as they go through the conjecture, produce the same *iterate* (*g* = *q*). This would require, from our earlier discussion, *p* to retain the binary form of *f* shifted to the left by two digits. However, we also realize that as the *3x+1* conjecture continues from *f*, it changes the binary form of *f* in the results that follow and <u>fails</u> to retain the original binary contents as it goes through the *Collatz trajectory* to reach *p*: *i.e.*, *p* will not have the binary structure of *f* and, as a result, $p \neq 4*f+1$, invalidating our assumption that *g* = *q*. (*Q. E. D.*)

**4.1. Generalization of 4x+1 Expression.** Now, let us consider a sequential version of Equation (4) with an intention to generalize Equation (1) that is applicable for *Terminal* integers, to other *odd* integers. Let $C_0$ and *y* be *odd* integers such that $y = (3C_0 + 1)/2^\alpha$. Also, let $C_k$ for *k* = 0, 1, 2, 3,… be such that $C_{k+1} = 4C_k + 1$. We present below an expression that lets all the $C_k$s, as they go through the Collatz conjecture, end up with the same *odd* integer $y = (3 C_k + 1)/2^\alpha$, each $C_k$ using an appropriate value for $\alpha$. This turns out to be the generalized version of Equation (1) that resembles Equation (4), (and here for *k* > 0):

$$C_k = 4^k C_0 + \sum_{i=0}^{k-1}(4^i) = 4^k C_0 + T_{k-1} \qquad (6)$$

Even though we mainly consider *odd* integers in this study, $C_0$ could be an *even* integer in Equation (6) as long as the application of Collatz conjecture is considered for $C_k$ with k > 0. Hence, the expression for $C_k$ in Equation (6) can be used for generating sets of *odd* integers with the *4x+1* characteristic. In the following, we use Equation (6) to generate all the sets of *odd* integers with *4x+1* characteristic so that each set's *Collatz iterate* is an *Intermediary* integer.

## 5. Intermediary Integers and their Predecessor sets

The main focus here is on the two types of *Intermediary* integers (6m+1 and 6m+5). In Section 1, we noted that the *Terminal* integers (1, 5, 21, 85, 341, …) that are linked by 4x+1 produce 1 as they go through the Collatz conjecture. Also, in Section 3, we noted that the *Pre-Terminal* integers of 5 consist of the following: 3, 13, 53, 213, … and these are also linked by 4x+1 with their Collatz *iterate* as 5. Here, we use these *Terminal* integers and *Pre-Terminal* integers (of 5) to start with and adapt Equation (6) suitably to find the *Collatz iterates* for all the *odd* integers.

The *iterate* for every *Terminal* integer is 1, which is also the first integer in the 6*m*+1 type *iterates* (*Intermediary* integers). Also, the *iterate* for the *Pre-Terminal* integers 3, 13, 53, … is 5, which is also the first integer in the 6*m*+5 type *iterates* (*Intermediary* integers). Related to the 4*x*+1 expression, we made an *Observation* in sub-section 3.2 that states: *relatively more odd integers, possibly 2/3$^{rd}$ of them, generate iterates of the type 6m+5 and the remaining 1/3$^{rd}$ only generate the iterates of the type 6m+1*. Our objective is to generate – with suitable $C_0$ in Equation 6 – all the remaining sets of *odd* integers (that are related by 4*x*+1) whose *Collatz iterates* are the remaining integers of both 6*m*+1 and 6*m*+5 (7, 13, 19, … and 11, 17, 23, …). To reflect the ratio of *odd* integers involved in producing the *iterates* of the two types, we could let $C_0$ = n for 6*m*+5 and $C_0$ = 2n for 6*m*+1, where n > 0. This could be one way of ensuring that 6*m*+5 will be the *iterates* for double the number of *odd* integers compared with the number of *odd* integers that 6*m*+1 will be the *iterates* for.

For illustration, in Table 4A, we start with just the first five *Terminal* integers 1, 5, 21, 85, and 341 (related by 4*x*+1) that are listed in the top row (n = 0) along with the first 6*m*+1 integer *i.e.,* 1 (for *m* = 0) in the last column. Also, in Table-4B, we start with only the first six *Pre-Terminal* integers of 5, which are: 3, 13, 53, 213, 853, and 3413 (related by 4*x*+1) that are listed in the top row (n = 0) along with the first 6*m*+5 integer *i.e.,* 5 (for *m* = 0) in the last column. In Table 4A, we create a sequence of *odd* integers for each of these five top row *Terminal* integers, listed under their respective columns, using the expression in Equation (6) with $C_0 = 2n$: $T_{k-1} + 4^k (2n)$; here *n* is the row (1, 2, 3,…) and *k* is the column (1, 2, 3, 4, 5). Since all the five *odd* integers in each row are related by *4x+1*, they all produce, using the respective *α* values shown in the Title row, the same *Collatz iterate* (listed in the last column of Table 4A). In Table 4B, we start with the first six *Pre-Terminal* integers for 5 and follow the same procedure except with $P_k$ in place of $T_k$ and $C_0 = n$ in Equation (6): $P_k + 4^k n$, where *n* is the row (1, 2, 3, …) and *k* is the column (1, 2, 3, 4, 5, 6). The six integers in each row produce, using the respective *α* values shown in the Title row, the same *Collatz iterate* (shown in the last column of Table 4B). We show only a few columns and rows in the Tables 4A & 4B; however *n* and *k* in the two expressions take all the positive integer values: $T_{k-1} + 4^k (2n)$ and $P_k + 4^k n$.

**5.1. Main Observations.** We list observations from the "*Extended*" Tables 4A and 4B.

*Observation 1:* The last column in Table 4A lists the *iterates* that are the *Intermediary* integers of the form *6m+1* for all *m* = 0, 1, 2, 3, …; and the last column in Table 4B lists the *iterates* that are the *Intermediary* integers of the form *6m+5* for all *m* = 0, 1, 2, 3, …. Hence each of these two last columns (listing *6m+ 1* and *6m+5*) have one-third of all the *odd* integers.

*Observation-2*: Examining the columns other than the two *iterate* ones (discussed above), we note that **all** the *odd* integers (including the *Starter* and the *Intermediary* integers) are present in the two *Extended* Tables 4A and 4B (**extended in both directions – *down* and *right*).** Thus **all** the *odd* integers collectively produce **all** the *Intermediary* integers as their *Collatz iterates*. Most importantly, there is **no duplication of any *odd* integer** in the two Tables.

*Observation-3:* Half of the *odd* integers are of the form *4m+3* ( m = 0, 1, 2, 3,…) and they are all under the second column of Table 4B, while *4m+1*, the remaining half, are spread out in all the remaining *non-iterate* columns of both the Tables 4A and 4B. We note that <u>one-third</u> of all the *odd* integers are in Table 4A: *i.e.,* adding all the fractions of *odd* integers in the various columns of Table 4A, we get $\sum_{i=1}^{\infty}(1/4)^i = 1/3$. As a result, <u>two-thirds</u> of all the *odd* integers (*i.e.*, all the 4*m*+3 and the remaining 4*m*+1) are in Table 4B.

.  11

| α in 2^α  Integer n | 2  1+8*n | 4  5+32*n | 6  21+128*n | 8  85+512*n | 10  341+2048*n | (3x+1)/2^α  For all Columns |
|---|---|---|---|---|---|---|
| 0 | 1 | 5 | 21 | 85 | 341 | **1** |
| 1 | 9 | 37 | 149 | 597 | 2389 | **7** |
| 2 | 17 | 69 | 277 | 1109 | 4437 | **13** |
| 3 | 25 | 101 | 405 | 1621 | 6485 | **19** |
| 4 | 33 | 133 | 533 | 2133 | 8533 | **25** |
| 5 | 41 | 165 | 661 | 2645 | 10581 | **31** |
| 6 | 49 | 197 | 789 | 3157 | 12629 | **37** |
| 7 | 57 | 229 | 917 | 3669 | 14677 | **43** |
| 8 | 65 | 261 | 1045 | 4181 | 16725 | **49** |
| 9 | 73 | 293 | 1173 | 4693 | 18773 | **55** |
| 10 | 81 | 325 | 1301 | 5205 | 20821 | **61** |
| 11 | 89 | 357 | 1429 | 5717 | 22869 | **67** |
| 12 | 97 | 389 | 1557 | 6229 | 24917 | **73** |
| 13 | 105 | 421 | 1685 | 6741 | 26965 | **79** |
| 14 | 113 | 453 | 1813 | 7253 | 29013 | **85** |
| 15 | 121 | 485 | 1941 | 7765 | 31061 | **91** |
| 16 | 129 | 517 | 2069 | 8277 | 33109 | **97** |
| 17 | 137 | 549 | 2197 | 8789 | 35157 | **103** |
| 18 | 145 | 581 | 2325 | 9301 | 37205 | **109** |
| 19 | 153 | 613 | 2453 | 9813 | 39253 | **115** |
| 20 | 161 | 645 | 2581 | 10325 | 41301 | **121** |
| 21 | 169 | 677 | 2709 | 10837 | 43349 | **127** |
| 22 | 177 | 709 | 2837 | 11349 | 45397 | **133** |
| 23 | 185 | 741 | 2965 | 11861 | 47445 | **139** |
| 24 | 193 | 773 | 3093 | 12373 | 49493 | **145** |
| 25 | 201 | 805 | 3221 | 12885 | 51541 | **151** |
| 26 | 209 | 837 | 3349 | 13397 | 53589 | **157** |
| 27 | 217 | 869 | 3477 | 13909 | 55637 | **163** |
| 28 | 225 | 901 | 3605 | 14421 | 57685 | **169** |
| 29 | 233 | 933 | 3733 | 14933 | 59733 | **175** |
| 30 | 241 | 965 | 3861 | 15445 | 61781 | **181** |
| 31 | 249 | 997 | 3989 | 15957 | 63829 | **187** |
| 32 | 257 | 1029 | 4117 | 16469 | 65877 | **193** |
| 33 | 265 | 1061 | 4245 | 16981 | 67925 | **199** |
| 34 | 273 | 1093 | 4373 | 17493 | 69973 | **205** |
| 35 | 281 | 1125 | 4501 | 18005 | 72021 | **211** |

TABLE 4-A.  Influence of 4x+1 in Collatz Conjecture Starting with Terminal Integers

| α in 2^α | 1 | 3 | 5 | 7 | 9 | 11 | (3x+1)/2^α |
|---|---|---|---|---|---|---|---|
| Integer n | 3+4*n | 13+16*n | 53+64*n | 213+256*n | 853+1024n | 3413+4096*n | for all columns |
| 0 | 3 | 13 | 53 | 213 | 853 | 3413 | **5** |
| 1 | 7 | 29 | 117 | 469 | 1877 | 7509 | **11** |
| 2 | 11 | 45 | 181 | 725 | 2901 | 11605 | **17** |
| 3 | 15 | 61 | 245 | 981 | 3925 | 15701 | **23** |
| 4 | 19 | 77 | 309 | 1237 | 4949 | 19797 | **29** |
| 5 | 23 | 93 | 373 | 1493 | 5973 | 23893 | **35** |
| 6 | 27 | 109 | 437 | 1749 | 6997 | 27989 | **41** |
| 7 | 31 | 125 | 501 | 2005 | 8021 | 32085 | **47** |
| 8 | 35 | 141 | 565 | 2261 | 9045 | 36181 | **53** |
| 9 | 39 | 157 | 629 | 2517 | 10069 | 40277 | **59** |
| 10 | 43 | 173 | 693 | 2773 | 11093 | 44373 | **65** |
| 11 | 47 | 189 | 757 | 3029 | 12117 | 48469 | **71** |
| 12 | 51 | 205 | 821 | 3285 | 13141 | 52565 | **77** |
| 13 | 55 | 221 | 885 | 3541 | 14165 | 56661 | **83** |
| 14 | 59 | 237 | 949 | 3797 | 15189 | 60757 | **89** |
| 15 | 63 | 253 | 1013 | 4053 | 16213 | 64853 | **95** |
| 16 | 67 | 269 | 1077 | 4309 | 17237 | 68949 | **101** |
| 17 | 71 | 285 | 1141 | 4565 | 18261 | 73045 | **107** |
| 18 | 75 | 301 | 1205 | 4821 | 19285 | 77141 | **113** |
| 19 | 79 | 317 | 1269 | 5077 | 20309 | 81237 | **119** |
| 20 | 83 | 333 | 1333 | 5333 | 21333 | 85333 | **125** |
| 21 | 87 | 349 | 1397 | 5589 | 22357 | 89429 | **131** |
| 22 | 91 | 365 | 1461 | 5845 | 23381 | 93525 | **137** |
| 23 | 95 | 381 | 1525 | 6101 | 24405 | 97621 | **143** |
| 24 | 99 | 397 | 1589 | 6357 | 25429 | 101717 | **149** |
| 25 | 103 | 413 | 1653 | 6613 | 26453 | 105813 | **155** |
| 26 | 107 | 429 | 1717 | 6869 | 27477 | 109909 | **161** |
| 27 | 111 | 445 | 1781 | 7125 | 28501 | 114005 | **167** |
| 28 | 115 | 461 | 1845 | 7381 | 29525 | 118101 | **173** |
| 29 | 119 | 477 | 1909 | 7637 | 30549 | 122197 | **179** |
| 30 | 123 | 493 | 1973 | 7893 | 31573 | 126293 | **185** |
| 31 | 127 | 509 | 2037 | 8149 | 32597 | 130389 | **191** |
| 32 | 131 | 525 | 2101 | 8405 | 33621 | 134485 | **197** |
| 33 | 135 | 541 | 2165 | 8661 | 34645 | 138581 | **203** |
| 34 | 139 | 557 | 2229 | 8917 | 35669 | 142677 | **209** |
| 35 | 143 | 573 | 2293 | 9173 | 36693 | 146773 | **215** |

TABLE 4-B.  Influence of 4*x*+1 in Collatz Conjecture Starting with Pre-Terminal Integers for 5

*Observation -4*: We also note that the *Starter* integers are evenly spread out in the two Tables 4A and 4B: *i.e.*, every third integer in all the columns (except the two *iterate* columns) and in all the rows is an *odd* multiple of 3. Since the *Starter* integers are evenly spread out, this explains why they use different values of α as they go through the Collatz conjecture process, as noted under *Observation (ii)* in Section 2. Here we clearly see how uniformly the *Starter* integers mix, in the Collatz conjecture process, with the *Intermediary* integers. It appears that the pending issue in [3] is to do with proving the kind of *mixing* between the two different categories of *odd* integers (*i.e., Starter* and *Intermediary*). We find this mixing to be uniform.

[Note: We now look at the significance of Tao's explanation for the "irregularities" in the behavior of these 6*m*+1 and 6*m*+5 *odd* integers (i.e., $\mathrm{Syr}(N)$ given below). Tao's statement from the first URL-link listed under [3] (that reflects *Observation-3*): *When viewed $3$-adically, we soon see that iterations of the Syracuse map become somewhat irregular. Most obviously, $\mathrm{Syr}(N)$ is never divisible by $3$. A little less obviously, $\mathrm{Syr}(N)$ is twice as likely to equal $2$ mod $3$ as it is to equal $1$ mod $3$. This is because for a randomly chosen odd $N$, the number of times $a$ that $2$ divides $3N+1$ can be seen to have a <u>geometric distribution</u> of mean $2$ – it equals any given value $a \in N+1$ with probability $2^{-a}$. Such a geometric random variable is twice as likely to be odd as to be even, which is what gives the above irregularity.*]

*Observation -5* In each row of the two Tables, there is a set of *odd* integers (related by *4x+1*) whose *Collatz iterate* – *i.e.*, (3x+1)/$2^\alpha$ using different values for $\alpha$ – match with the *Intermediary* integer in the last column of that row. Thus each row displays the entire *predecessor set* of the corresponding *Intermediary* integer in the last column, addressing the issues discussed on pages 62-63 in [2] regarding the size of *predecessor* set.

*Observation-6*: The Tables 4A and 4B provide the actual evidence for the *heuristic* probabilities discussed on page 34 in [2] (and also the URL listed under [2]) that *the expected growth in size between two consecutive odd integers in such a trajectory is the multiplicative factor ¾ < 1; i.e., on the average, the iterates in a trajectory tend to shrink in size.* We see the following evidence in the two Tables. In the second column of Table-4B, all the alternate *odd* integers (*i.e.,* ½ of all the *odd* integers that are of the form 4*m*+3, for all *m* = 0, 1, 2, 3, …) are listed and they all use $\alpha$ = 1 in (3*x*+1)/$2^\alpha$. Hence each of these integers in this column produces a conjecture result that is larger (by about 1 ½ times) than itself. (In Section 7, we conduct an in depth analysis of the 4*m*+ 3 type integers.) In the two Tables, all the *odd* integers in all the other columns, that are of type 4*m*+1, use $\alpha \geq 2$ in (3*x*+1)/$2^\alpha$ to produce the results shown in their respective last column. Specifically, looking at all the columns in the two Tables, we note that, as the value of $\alpha$ increases by 1, the number of *odd* integers in their respective column reduces by half (*i.e.*, higher the value of $\alpha$, proportionately fewer the number of *odd* integers). We will have more discussion on these observations in Section 8. Now, we show, in the next sub-section, that the *Starter* integers are like the *leaves* on the *Collatz tree*

**5.2. Constructing Collatz Tree and Conjecture Trajectory.** We transform the *Collatz Graphs* (Figures 1A & 1B) into a single *Collatz Tree,* shown in Figure 2, with its root (*i.e.*, *End* integer 1) above, the branches (*i.e.*, *Intermediary* integers) below and the leaves *(i.e., Starter* integers) spread out everywhere. Here, we also use some results from the *Extended* Tables 4A and 4B and also the 4*x*+1 expression to include additional *odd* integers so that it is symmetrical. Also in Figure 3, we show a *Collatz Tree* with just the *iterates* (*i.e., the Intermediary* integers). In Section 6, we discuss a method for systematically generating <u>all</u> the *odd* integers for the <u>orderly placement</u> in different lower layers of the *Collatz Tree* with 1 at the top.

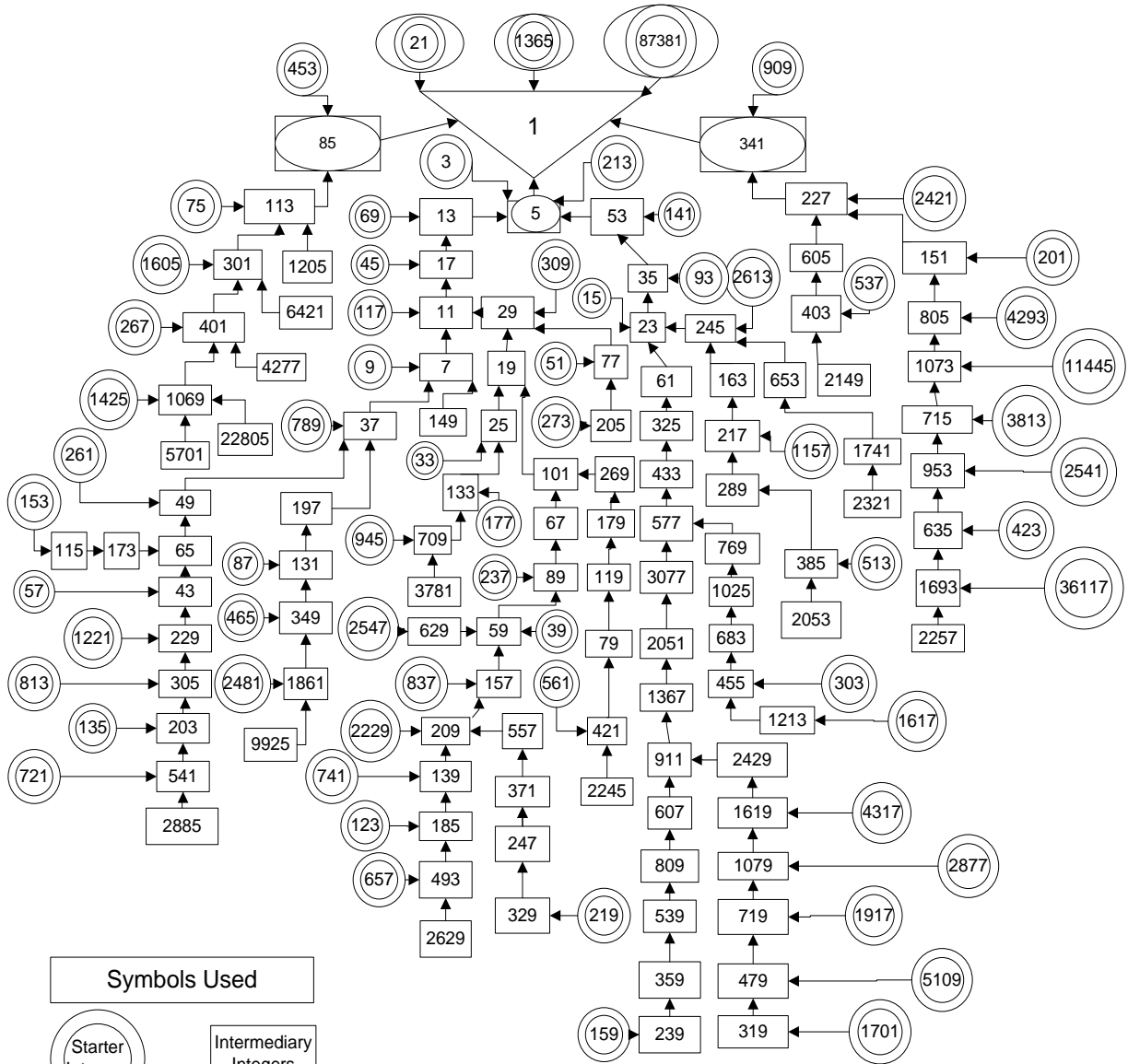

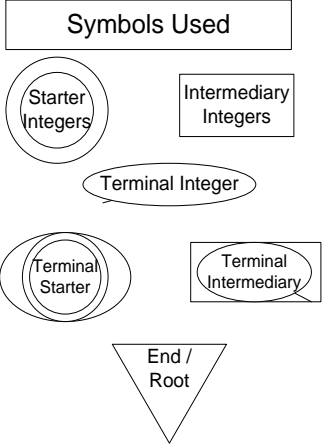

Symbols Used

Figure 2: Example 3x+1 Tree

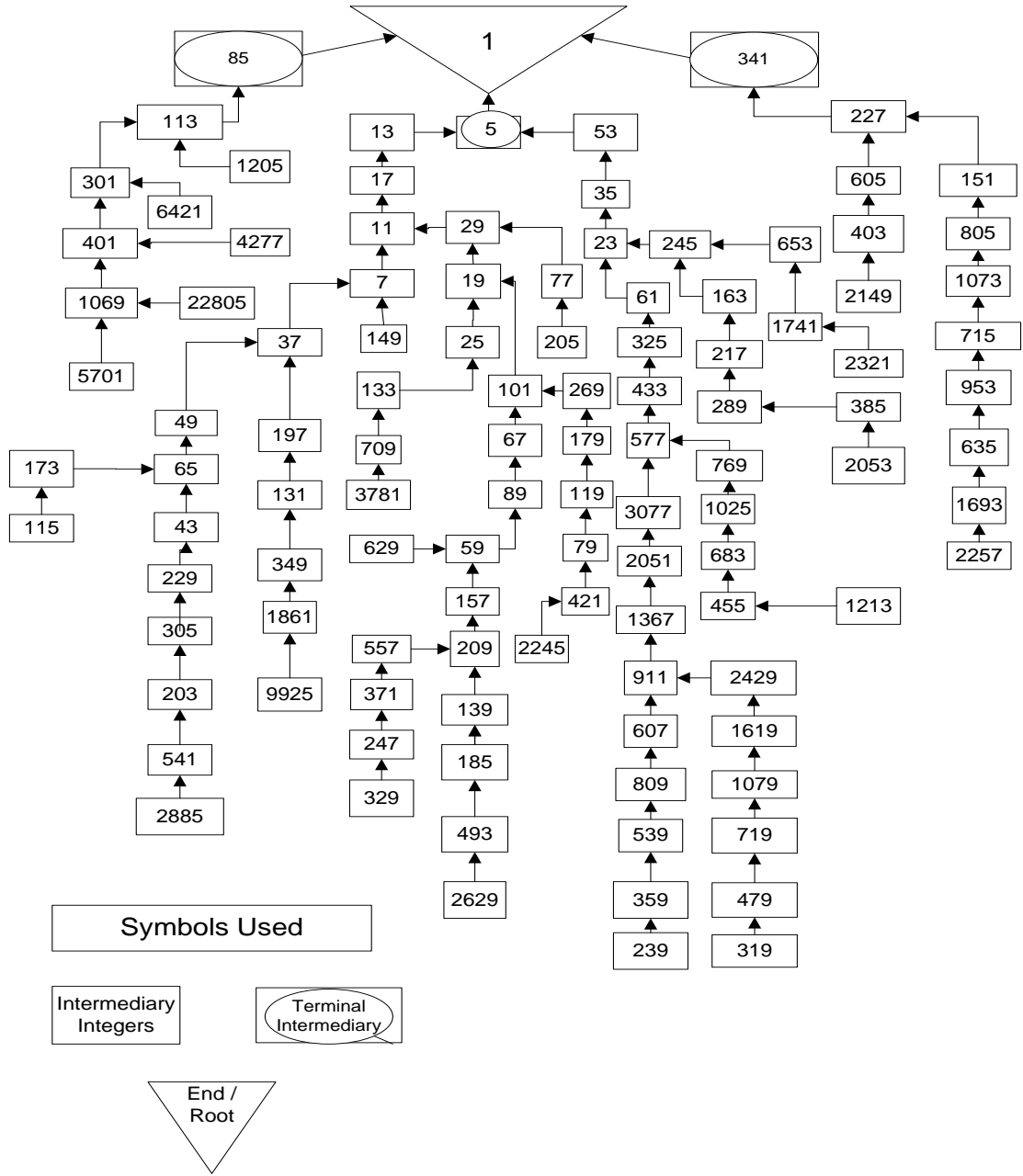

Figure 3: Example 3x+1 Tree With Intermediary Integers Only

Here is a simple logic to build the Collatz *trajectory* for a given *odd* integer $x > 1$. Using the rows of *odd* integers along with their conjecture *iterate* that are in Tables 4A and 4B, it is possible to construct the *Conjecture Trajectory* for any given *odd* integer in a simpler way. Given an *odd* integer $x$ (other than 1) we first find where $x$ is positioned in the Tables 4A and 4B: *i.e.,* (i) in which of the two Tables and (ii) in which row in that Table so that we can get its *Collatz iterate* (say $y$) from the last column of the corresponding Table.

We begin by checking if it is an integer of the form $4m+3$ (*i.e.,* half of the *odd* integers) and these are <u>all</u> in the first column of Table 4B. Among the $4m+1$ type integers, half of them are in the first column of Table 4A and rest are distributed. Hence, as the second possibility, we check if it is in the first column of Table 4A. If the $4m+1$ integer is in some other column, we look for the integer in the previous column of the row where $x$ is and make that integer as the new $x$ and repeat the process until x is in one of the first columns in the two Tables. First we print the given x and as we go through the process, we only print all those *iterates* (values of y) when the x is in one of the first columns. Finally when x becomes 1, we print 1 and end the process. The simple logical steps are presented below:

Print x;
Begin: print '→' if [x = 1]
   then print x; Stop
   else if [x(mod)4 = 3]
      then [ m = ⌊x/4⌋; y = 6*m+5; print y; x = y; Goto Begin;] ;// first column of Table 4B
      else if [x(mod)8 = 1]
         then [ m = ⌊x/8⌋; y = 6*m+1; print y; x = y; Goto Begin;] //first column of 4A
         else x = (x–1)/4; Goto Begin; // in some column other than the first in 4A or 4B

## 6. Mechanism for Constructing a Collatz Tree

In the *Collatz Tree* shown in Figure 3, the branches are the *Intermediary* integers which constitute the complete set of solutions $(3x+1)/2^\alpha$ (samples presented in Tables 4A and 4B) where $x$ represents every *odd* integer and these *branches* lead up to the *root* at the top (*i.e.*, *End* integer 1). We skip the *Starter* integers since, as they go through the conjecture process, they get attached as *leaves* (by not being *iterates* for any *odd* integers). Hence they are ignored when creating the next lower layer.

The *extended* Tables 4A and 4B together list all the *Intermediary* integers as the *Collatz conjecture iterates* for sets of connected *odd* integers shown in their respective rows. Here we provide a procedure that could be mechanized for drawing a *Collatz Tree,* starting with 1 at the top (Layer-0). The next layer, Layer-1, consists of all the *Terminal* integers (5, 21, 85, 341, 1365, 5461, …) whose Collatz conjecture result is 1 and hence they are all linked upwards to 1. The layer below this consists of the sets of *Pre-terminal* integers for each of these *Terminal* integers (except for the *Starter* integers like 21, 1365, ). So, Layer-2 consists of the *Pre-terminal* integers under each of the *non-Starter Terminal* integers. For <u>5</u> they are: (3), 13, 53, (213), 853, 3413,… ; for <u>85</u>: 113, (453), 1813, 7253,…; for <u>341</u> are: 227,(909), 3637, 14549,…; and so on (those in parenthesis are *Starter* integers). We see that the *Collatz Tree* spreads out as the process goes down the layers. The initial set up is given here with $i = 1, 2,, 3,$ indicating the sequence: L0 = 1; L1($i$) = 5, 85, 341, 5461,…; L2(1, $i$) = 13, 53, 853,…; L2(2, $i$) = 113, 1813, 7253,…; L2(3,$i$) = 227, 3637, 14549,…; and so on. We note that the results in the first row of the two extended Tables 4A and 4B appear in layers 1 and 2 respectively. For the remaining lower layers we use, given a *Collatz iterate*, a simple mechanism for finding the corresponding lower layer *odd* integers. This mechanism consists of a few steps and we start with these results in the third layer: L2 and use that as the new *y* to continue the process. The mechanism that could automatically get such next *y* is described next.

**6.1 Collatz Tree Construction Mechanism.** Let *y* be an *odd* integers from Layer-2, say *y* = L(*n, j*). We first find where *y* is positioned in the Tables 4A and 4B: *i.e.,* (i) in which of the two result-columns and (ii) in which row of that column so that we can get the set of *odd* integers whose *Collatz iterate* is *y*. We place all these *odd* integers in the layer just under *y* and continue the process by picking the first *odd* integer in the list (that is not a *Starter* integer) and make that as the new *y* and repeat the process. The first entries in the first columns of the two Tables are different (1 for Table 4A and 3 for Table B). Since we need these integers for later use, we also use them to identify the Tables. Also, the row can be identified by computing the integer part of y/6. This logic is presented below:

If [y(mod)6 = 1] then s=1 else s = 3
        // s= 1 top entry in the first column. in Table 4A & s= 3 in the first column in Table 4B
        // We note, 13 is in Table 4A; 53 is in Table 4B which are used as examples below
r = ⌊y/6⌋      // Here *r* represents the row
        // We note that 13 is in row 2 (of Table 4A) and 53 is in row 8 (of Table B)
        // At this stage we know the Table and the row where the Collatz conjecture result y is
        // If the current y is in Table 4A we use the following;
        // also make the first *non-Starter* integer as the next y
If s = 1 then {z = 8*r + 1;  y = z;  if [y (mod) 3] = 0 then [y = 4*z + 1]  // z is in Table 4A
        // Now, create an array e(n) to store the first *n odd* integers that are in row *r* of Table 4A
        e(1) = z; do for all *i* = 1 to n        [e(*i*+1) = 4* e(*i*) + 1]}
        // If the current y is in Table 4B, we use the following;
        // also make the first *non-Starter* integer as the next y
If s = 3 then {z = 4*r + 3; y = z; if [y (mod) 3] = 0 then [y = 4*z + 1] // z is in Table 4B
        // Now, create an array e(*n*) to store the first *n odd* integers in row *r* of Table 4B
        e(1) = z; do for all *i* = 1 to *n*        [e(*i*+1) = 4* e(*i*) + 1]}  //
        // Store or use these results and move on to the next layer with the new value of *y*.

In the above, assuming we started with 13, we note z = 17 and hence the new y = 17; e(n) will contain 17, (69), 277,…; *z* and e(n) will go in the layer below 13 in the *Collatz tree*. Similarly with 53, z = 35, the new y = 35; e(n) will contain 35, (141), 565,…. As a result, Layer-3 will have several *segments* linking up to 13: we have L3(1,1,*i*) = 17, (69), 277, …; and also segments linking up to 53: L3(1,2,*i*) = 35, (141), 565,…Obviously , for drawing the *Collatz Tree*, the above may not be the most efficient approach and certainly it is not comprehensive. However, our objective here is just to indicate how the results in the Tables 4A and 4B can be interpreted for generating the <u>entire</u> *Collatz Tree* leading up to the *root* integer 1 at the top.

## 7. Collatz Trajectory Characteristics of 4m+3 Integers

Each of the *4m+3* integers uses $\alpha = 1$ at least once as it starts going through the Collatz conjecture process, producing an *iterate* that is larger (by about 1 ½ times ) than itself. Some of such *Collatz iterate*s may also use $\alpha = 1$ in their Collatz conjecture process and this pattern could continue. For instance, we see in the *Collatz trajectory* for 63 (which is a *4m+3* integer) presented in Table 2A: **63** → 95 → 143 → 215 → 323 → 485 → 91 → … → 1. There are five increases before a reduction since 63, 95, 143, 215, and 323 are all integers of type *4m+3* (and use $\alpha = 1$) while 485, an integer of type *4m+1* using $\alpha = 4$, produces a lower *iterate*.

As a result, there is a need to classify the *4m+3* integers listed under the second column in Table 4B based on their length of continuously increasing *Collatz trajectory* (for which we use the term: *Alpha1 Trajectory Length*, since we are considering all the integers using $\alpha = 1$ in the *trajectory*). For instance, the *Alpha1 Trajectory Length* for 63 is 5 as each of these five integers use $\alpha = 1$: **63** → 95 → 143 → 215 → 323. Next, we present a method for categorizing the entire *4m+3* integers based on their *Alpha1 Trajectory Lengths*.

**7.1 Alpha1 Trajectory Lengths of 4m+3 Integers.** Let $y_1 = (3x+1)/2$, where *x* and *y* are *odd* integers and *y* is *x*'s first *Collatz iterate* (using $\alpha = 1$). When *x* happens to have a higher *Alpha1Trajectory Length*, we will use, depending on its *trajectory* length, one of these equations: $y_2 = (3y_1+1)/2 = (9x+5)/4$; or $y_3 = (27x+19)/8$; or $y_4 = (81x+65)/16$; or $y_5 = (243x+211)/32$; and so on. Let *h* denote the *Alpha1 Trajectory Length*. First, find a set of *odd* integers with *h* =1 using $y_1 = (3x+1)/2$. To ensure that no *odd* integers with h > 1 are included, one approach would be to consider integers of *4m+1* whose Collatz conjecture, we know, will use $\alpha \geq 2$ (seen in Figures 1A and 1B). Hence, if $y_1 = 4z + 1$ where z is an *odd* integer, we note that $(3y_1+1)/2 = [3(4z + 1) +1]/2 = (12z +4)/2$, indicating that the conjecture result for $y_1$ will be using $\alpha \geq 2$ and hence this conjecture result cannot belong to the category of *odd* integers with $\alpha = 1$ (ensuring that *x* belongs to the category of *h* = 1 only and not to any higher values of *h*)

The above being a desired result for what we are looking for, the *odd* integers with *h* = 1 could belong to the group where $y_1 = (3x+1)/2 = 4z+1$. So, $z = (3x-1)/8$, which can be used to find all the *odd* integers *x* with *h* =1 while ensuring that *z* remains an *odd* integer. Such *odd* integers of *x* are: 3, 11, 19, …; *i.e.*, *x* = 3+8(*n* – 1), where *n* = 1, 2, 3, ... Following this approach, it is possible to get the set of all *odd* integers with *h* = 2 only by using $y_2 = (9x+5)/4 = 4z+1$. The *odd* integers for *h* = 2 are chosen so that the integer $z = (9x+1)/16$ and they are: 7, 23, 39, 55, …; *i.e.*, *x* = 7 +16(*n* –1), where *n* = 1, 2, 3,…; further, for *h* = 3, z = (27*x*+11)/32 and the respective *odd* integers are: 15, 47, 79,…; *i.e.*, *x* = 15 +32(*n* –1), where *n* = 1, 2, 3,…; and this process could be continued, using the appropriate expression for any $y_h$, to find the *odd* integers with higher values of *h*.

We note that the first integer in each of the above three series show a pattern: $2^{h+1} - 1$. This is not a new pattern in the study of 3*x*+1, since we see reported results on page 34 in [2] from the earlier studies on *odd* integers such as $2^{50}-1$ and $2^{500}-1$. Since the current study is about *Alpha1 Trajectory Lengths* we also provide, for different values of *h,* the set of integers that follow each of the first integers in $2^{h+1} - 1$. The list of *odd* integers for the increasing number of *Aplha1 Trajectory Lengths* (*i.e.*, *h* = 1, 2, 3, …) are presented in Table-5, and this Table, just as Tables 4A and 4B, could be **extended both to the *right* and *down*** so as to include **all** the integers of the type 4*m*+3. The extended entries in Table-5 (*i.e.*, for all positive integer values of *h and n*) is generated from the following comprehensive and generalized expression that captures the earlier observations, where *n* ($\geq 1$) indicates the row:

$$x_{h,n} = 2^{h+1}(2n - 1) - 1 \qquad (7)$$

*Main Observation:* In the above, the focus was on generating comprehensive results concerning *Alpha1 Trajectory Lengths* of 4*m*+3. In Table-5, we note that the first column (*with Aplha1 trajectory length* =1) has half of these 4*m*+3 integers and, as the length increases, the number of entries in each of the subsequent columns reduces by half (a similar pattern we saw in Tables 4A and 4B), indicating that these three Tables contain comprehensive information regarding the actual behavior of the Collatz conjecture.

| α=1 Chain Length | 1 | 2 | 3 | 4 | 5 | 6 | 7 | 8 | 9 | 10 |
|---|---|---|---|---|---|---|---|---|---|---|
| List n= | $x_{1,n}$ | $x_{2,n}$ | $x_{3,n}$ | $x_{4,n}$ | $x_{5,n}$ | $x_{6,n}$ | $x_{7,n}$ | $x_{8,n}$ | $x_{9,n}$ | $x_{10,n}$ |
| 1 | 3 | 7 | 15 | 31 | 63 | 127 | 255 | 511 | 1023 | 2047 |
| 2 | 11 | 23 | 47 | 95 | 191 | 383 | 767 | 1535 | 3071 | 6143 |
| 3 | 19 | 39 | 79 | 159 | 319 | 639 | 1279 | 2559 | 5119 | 10239 |
| 4 | 27 | 55 | 111 | 223 | 447 | 895 | 1791 | 3583 | 7167 | 14335 |
| 5 | 35 | 71 | 143 | 287 | 575 | 1151 | 2303 | 4607 | 9215 | 18431 |
| 6 | 43 | 87 | 175 | 351 | 703 | 1407 | 2815 | 5631 | 11263 | 22527 |
| 7 | 51 | 103 | 207 | 415 | 831 | 1663 | 3327 | 6655 | 13311 | 26623 |
| 8 | 59 | 119 | 239 | 479 | 959 | 1919 | 3839 | 7679 | 15359 | 30719 |
| 9 | 67 | 135 | 271 | 543 | 1087 | 2175 | 4351 | 8703 | 17407 | 34815 |
| 10 | 75 | 151 | 303 | 607 | 1215 | 2431 | 4863 | 9727 | 19455 | 38911 |
| 11 | 83 | 167 | 335 | 671 | 1343 | 2687 | 5375 | 10751 | 21503 | 43007 |
| 12 | 91 | 183 | 367 | 735 | 1471 | 2943 | 5887 | 11775 | 23551 | 47103 |
| 13 | 99 | 199 | 399 | 799 | 1599 | 3199 | 6399 | 12799 | 25599 | 51199 |
| 14 | 107 | 215 | 431 | 863 | 1727 | 3455 | 6911 | 13823 | 27647 | 55295 |
| 15 | 115 | 231 | 463 | 927 | 1855 | 3711 | 7423 | 14847 | 29695 | 59391 |
| 16 | 123 | 247 | 495 | 991 | 1983 | 3967 | 7935 | 15871 | 31743 | 63487 |
| 17 | 131 | 263 | 527 | 1055 | 2111 | 4223 | 8447 | 16895 | 33791 | 67583 |
| 18 | 139 | 279 | 559 | 1119 | 2239 | 4479 | 8959 | 17919 | 35839 | 71679 |
| 19 | 147 | 295 | 591 | 1183 | 2367 | 4735 | 9471 | 18943 | 37887 | 75775 |
| 20 | 155 | 311 | 623 | 1247 | 2495 | 4991 | 9983 | 19967 | 39935 | 79871 |
| 21 | 163 | 327 | 655 | 1311 | 2623 | 5247 | 10495 | 20991 | 41983 | 83967 |
| 22 | 171 | 343 | 687 | 1375 | 2751 | 5503 | 11007 | 22015 | 44031 | 88063 |
| 23 | 179 | 359 | 719 | 1439 | 2879 | 5759 | 11519 | 23039 | 46079 | 92159 |
| 24 | 187 | 375 | 751 | 1503 | 3007 | 6015 | 12031 | 24063 | 48127 | 96255 |
| 25 | 195 | 391 | 783 | 1567 | 3135 | 6271 | 12543 | 25087 | 50175 | 100351 |
| 26 | 203 | 407 | 815 | 1631 | 3263 | 6527 | 13055 | 26111 | 52223 | 104447 |
| 27 | 211 | 423 | 847 | 1695 | 3391 | 6783 | 13567 | 27135 | 54271 | 108543 |
| 28 | 219 | 439 | 879 | 1759 | 3519 | 7039 | 14079 | 28159 | 56319 | 112639 |
| 29 | 227 | 455 | 911 | 1823 | 3647 | 7295 | 14591 | 29183 | 58367 | 116735 |
| 30 | 235 | 471 | 943 | 1887 | 3775 | 7551 | 15103 | 30207 | 60415 | 120831 |
| 31 | 243 | 487 | 975 | 1951 | 3903 | 7807 | 15615 | 31231 | 62463 | 124927 |
| 32 | 251 | 503 | 1007 | 2015 | 4031 | 8063 | 16127 | 32255 | 64511 | 129023 |
| 33 | 259 | 519 | 1039 | 2079 | 4159 | 8319 | 16639 | 33279 | 66559 | 133119 |
| 34 | 267 | 535 | 1071 | 2143 | 4287 | 8575 | 17151 | 34303 | 68607 | 137215 |
| 35 | 275 | 551 | 1103 | 2207 | 4415 | 8831 | 17663 | 35327 | 70655 | 141311 |
| 36 | 283 | 567 | 1135 | 2271 | 4543 | 9087 | 18175 | 36351 | 72703 | 145407 |

TABLE 5  Integers Producing Collatz Results that Continuously Increase (First 10 $\alpha 1$ Chain Length)

We present, in the following Section, the evidence for *on-the- average* convergence of <u>all</u> the *odd* integers' *Collatz conjecture trajectories* using the results presented in Tables 4A, 4B and 5.

## 8. Validity of the Collatz Conjecture

Earlier in Section 5 under *Observation* 6, we used a quote from page 34 in [2]. By expanding the idea in that quote, we now present a Theorem that captures the essence of the Collatz conjecture. Here we consider the mixing properties of the *iterates* of the two types of *odd* integers: (i) 4m+3 whose *iterates* keep increasing and 4m+1 whose *iterates* keep decreasing.

THEOREM 8.1 The Collatz conjecture results of half of all the *odd* integers that are of the form 4*m*+3, on the average, increase by a factor of 3 and the results of the other half of *odd* integers that are of the form 4*m*+1, on the average, decrease by a factor 4 (*i.e.,* a multiplication factor of ¼). As a result, considering all the *odd* integers, there are no divergent *Collatz conjecture trajectories*.

*Proof*: All the 4*m*+3 *odd* integers are listed under the column $\alpha = 1$ in Table 4B. In Table 5, these (4*m*+3) *odd* integers are, based on their value of *Aplha1 Chain Lengths (h)*, distributed under different columns of *h* (*i.e., h* = 1, 2, 3, …)

In Table 5, half of the 4*m*+3 type integers are listed under the first column $h = 1$ and the number of integers under each of the subsequent columns reduces by half as the value of *h* increases by 1. The *Collatz iterate* for each of the *odd* integer under the first column increases roughly by 3/2 – ignoring the 1 in the expression (3x+1)/2. Similarly the *Collatz iterates* in the subsequent columns also increase by roughly $(3/2)^h$. The overall average increase is given by:

$$(1/2)(3/2) + (1/4)(3/2)^2 + (1/8).(3/2)^3 + \ldots = (3/4) + (3/4)^2 + (3/4)^3 \ldots$$
$$= \sum_{n=1}^{\infty}(3/4)^n = \mathbf{3} \qquad (8)$$

This could be interpreted as if *the Collatz iterates* of half of all the *odd* integers that are of the form 4*m*+3, on the average, increase by <u>three times</u> the value of the integer that was used to start with. Now let us examine all the remaining columns ($\alpha \geq 3$) in Table 4B that has one-third of the 4*m*+1 *odd* integers. The *Collatz iterate* for each of the *odd* integer under the column $\alpha = 3$ decreases roughly to (3/8) times its original value and there are 1/8[th] of all the *odd* integers here.. The average decrease of these one-third of the 4*m*+1 integers is given by:

$$(1/8)(3/8) + (1/32)(3/32) + (1/128)(3/128) + \ldots = (3/8^2) + (3/32^2) + (3/128^2) \ldots$$
$$= (3/2^2)\sum_{n=1}^{\infty}(1/4^2)^n \qquad (9)$$

We now consider the remaining two-thirds of the 4m+1 *odd* integers listed under all the columns in Table 4A with *even* values of $\alpha$. The *Collatz iterate* for each of the integer under the column $\alpha = 2$ decreases roughly to (3/4) times its original value and there are 1/4[th] of all the *odd* integers here  The average decrease of the remaining two-thirds of the 4*m*+1 integers is given by:

$$(1/4)(3/4) + (1/16)(3/16) + (1/64)(3/64) + \ldots = [3/(4)^2] + [3/(16)^2] + [3/(64)^2]\ldots$$
$$= 3\sum_{n=1}^{\infty}(1/4^2)^n \qquad (10)$$

Combining the results in Equations (9) and (10), we note that the *Collatz iterates* of the other half of *odd* integers that are of the form 4*m*+1, on the average, decrease and their overall multiplication factor is given by:

$$(3/2^2)\sum_{n=1}^{\infty}(1/4^2)^n + (3)\sum_{n=1}^{\infty}(1/4^2)^n = (15/4)\sum_{n=1}^{\infty}(1/4^2)^n$$
$$= (15/4)(1/15) = \mathbf{1/4} \qquad (11)$$

This could be interpreted as if the *Collatz iterates* of <u>half</u> of the *odd* integers of the form 4*m*+1, on the average, get smaller by a factor <u>4</u> of the integer that was used to start with.

In the above, we have the results for the *odd* integers of the type 4*m*+3 in Equation(8) that shows an average *increase by a factor of 3* and for those of the type 4*m*+1 in Equation(11) that shows an average *decrease by a factor of 4*. Both have equal number of *odd* integers in them. [Note: The *Observation*-4 under *Sub-Section 5.1* is relevant here: The *Starter* integers are evenly spread out in the two Tables 4A and 4B: *i.e.*, every third integer in all the columns (except the two *iterate* columns) and in all the rows is an *odd* multiple of 3. In view of this uniform distribution of *the Starter* integers, we can assume that they are uniformly excluded in all the Equations (8) to (11) without affecting the final results]. The overall result (3/4) matches with the theoretical (heuristic) result 0.75 on page 34 in [2]. Hence, on the average, the *odd* integer *iterates* in the *Collatz conjecture trajectory* shrink and hence there cannot be any divergent *Collatz conjecture trajectories*. *(Q. E. D.)*

## 9. Conclusion

First, we summarize the classification of *odd* integers used in this study on Collatz conjecture. *Terminal* integers are those whose Collatz conjecture result is 1, the *End* integer. We classified all the *odd* integers in two ways: (i) consisting of three equal categories: *odd* multiples of 3 (known as *Starter* integers); 6*m*+1, and 6*m*+5, with *m* =0, 1, 2, 3,… (the last two collectively known as *Intermediary* integers); and (ii) consisting of two equal categories: 4*m*+1 and 4*m*+3 (where *m* =0, 1, 2, 3,…). The first type of classification is used in [3] without these names. In Section 5, we clearly explain the *mixing* between these two types of integers (*Starter* and *Intermediary*) in *Observation 4* that pertains to all the *odd* integers in Tables 4A and 4B.

We started with the *Collatz trajectories* for a few *Starter* integers and, using these results, we presented a new type of *Collatz Graph*. From this graph, we identified the role that 4*x*+1 expression plays in the analysis of the Collatz conjecture. Further in Theorem 4.1, we proved that a *Collatz trajectory* cannot have duplicate *odd* integers in it (hence no looping as well). One-half of the *odd* integers of type $x = 4m+1$ use values of $\alpha \geq 2$ in the expression $(3x+1)/2^\alpha$ while the other half of the *odd* integers of type $x = (4m+3)$ just use $(3x+1)/2$ (*i.e.*, $\alpha = 1$).

The valuable role that the <u>generalized</u> version of 4*x*+1 expression plays in this study was demonstrated through the creation of Tables 4A and 4B (and also, to some extent, Table 5). The *Intermediary* integers of type *6m+1* (last column in Table 4A) are the *Collatz iterates* of <u>one third</u> of the *odd* integers and also those of type *6m+5* (last column in Table 4B) are the *Collatz iterates* of the remaining <u>two thirds</u> of all the *odd* integers. As can be seen, we have avoided the need for tools such as the *Syracuse random variables* in order to identify all the *iterates* that are listed in the last two columns of the Tables 4A and 4B. In Table 4B, the 4*m*+3 integers make half of all *odd* integers that are listed under the column $\alpha = 1$; and one-sixth of the *odd* integers in Table 4B are of the type 4*m*+1 that are listed under the remaining columns $\alpha \geq 3$. The remaining 4*m*+1 integers are presented in Table 4A shown in columns with *even* values of $\alpha$. We also looked into the behavior of these 4*m*+3 integers in terms of their *Alpha1 trajectory length* and presented these results in Table 5 in the increasing order of *Alpha1 trajectory lengths*.

From the proof of the Collatz conjecture perspective, Tables 4A and 4B (excluding the column for $\alpha = 1$) show how well, from Collatz conjecture perspective, the integer *(iterate)* values get reduced, while Table 5 shows how badly, again from Collatz conjecture perspective, the 4*m*+3 integer *(iterate)* values get increased. In Theorem-8.1, we proved that the Collatz conjecture *iterate* of all the *odd* integers considered together shrink, on the average, by a factor of 3/4 (matching with the *heuristic* result reported in [2]). The ideas used here are rather simple. However, we go into depths when analyzing the conjecture showing the mixing properties of the *iterates* of the two types of *odd* integers.

Since it is difficult to trace all the earlier contributions, we used [2] as our main source as it presents all the past information in a condensed manner. The recent contribution by Tao [3], where the conjecture has been *almost* proved, takes a highly formal approach. Even though one of the categorization of *odd* integers used here is similar to the one seen in [3], we have taken a structured approach to observe all the hidden intricacies of the conjecture. Here, all the *odd* integers are organized based on their conjecture *iterates*. We have described a simple logic involving expressions for locating integers in Tables 4A and 4B [and not using again the conventional $(3x+1)/2^\alpha$] for generating the Collatz *trajectory* for any given *odd* integer $x > 1$. Similarly, we have come up with a mechanism, based on the results in Tables 4A and 4B, to draw the complete *Conjecture Tree* starting with 1 at the top, the *Terminal* layer (consisting of *Terminal* Integers) below that, and the appropriate *Pre-terminal* layer below each of the *Terminal* integer. This process can be used to construct the *Collatz conjecture net (tree)* without resorting to further use of $(3x+1)/2^\alpha$ expression. In Theorem 8.1, considering all the *odd* integers and their conjecture *iterates* we prove that there cannot be any divergent Collatz conjecture *trajectories*.

The current practice used in solving the Collatz conjecture seems to be based on a verbatim interpretation, where we start with an integer and produce a trajectory to reach the *end* integer 1. However, after using the **4x+1 key** to open the *(3x+ 1)'s-box* and seeing what are all inside (such as samples displayed in Tables 4A, 4B and 5), we may have to come up with a (w)holistic approach. The current interpretation is like looking at strings (*trajectories*) of a big *net* instead of recognizing the *net* itself. The new Theorem statement, perhaps, could be:

*All odd integers like x and y that are related by $(3x+1)/2^\alpha = y$, where $\alpha$ is a suitable integer, go up through the Collatz net with y in a layer above x to reach 1.*

*Acknowledgement:* We thank Dr. Jeffrey C. Lagarias for pointing out that (i) the categorization of *odd* integers into two main types (similar to *Starter* and *Intermediary* integers used in this paper) has already been discussed by Dr. Tao in [3] and (ii) there is still a need to prove the kind of *mixing* between these two categories of *odd* integers.

*Kamayasamy (Ken) Surendran* retired from Southeast Missouri State University as Professor of Computer Science in 2013.  His previous academic assignments in Computer Science were with Rose-Hulman Institute of Technology; UNITEC Institute of Technology, New Zealand; and PSG College of Technology, India.  His industrial experiences in IT were with Indian Space Research Organization and Zambia Consolidated Copper Mines.  Surendran received B.E. (1965) in Electrical Engineering from University of Madras, India, M. Tech. (1967) in Electrical Engineering (Control Systems) from Indian Institute of Technology, Madras, India, and Ph. D. (1971) in Applied Analysis / Applied Mathematics from State University of New York at Stony Brook.  ( ksurendran@semo.edu or suren@linuxmail.org )

*Desrazu Krishna Babu* retired from Mobil Oil Corporation in 2000 after serving as Senior Research Engineer (R&D) for fifteen years.  Earlier he was with the Civil Engineering Department, Princeton University, as Research Engineering Staff; and with the Mathematics Department, City College, CUNY, as Assistant Professor. Babu received M.A. (1962) in Applied Mathematics from Osmania University, India, and Ph.D. (1971) in Applied Analysis /Applied Mathematics from State University of New York at Stony Brook. ( sithababu39@gmail.com )